\documentclass[A4j,11pt]{article}
\usepackage{amssymb,amsmath,amsthm,amscd,amsfonts}
\usepackage{amsmath,amssymb,amsthm,amscd,ascmac,mathrsfs,bm,type1cm,multirow,array,enumerate,mathrsfs,tabularx}
\usepackage{graphicx}
\usepackage{graphics}
\newcommand{\vu}{{\bm{u}}} %m変数
\newcommand{\vv}{{\bm{v}}} %ベクトル場
\newcommand{\vN}{\bm{N}} %法線ベクトル

\newcommand{\vY}{\bm{Y}}

\newcommand{\vV}{\bm{V}}

\begin{document}

\title{{\bf A construction of \\
curvature-adapted hypersurfaces \\
in the product of symmetric spaces}}
\author{{\bf Naoyuki Koike}}
\date{}
\maketitle

\begin{abstract}
In this paper, we give a construction of curvature-adapted hypersurfaces in the product $G_1/K_1\times G_2/K_2$ of (Riemannian) symmetric spaces 
$G_i/K_i$ ($i=1,2$).  By this construction, we obtain many examples of curvature-adapted hypersurfaces in $G_1/K_1\times G_2/K_2$.  
Also, we calculate the eigenvalues of the shape operator and the normal Jacobi operator of the curvature-adapted hypersurfaces obtained 
by this construction.  
\end{abstract}

%\subjclass{53D12, 53C35}

\vspace{0.5truecm}

\section{Introduction}
Let $M$ be a connected oriented hypersurface in a Riemannian manifold $(\widetilde M,\widetilde g)$ immersed by an immersion $f:M\hookrightarrow\widetilde M$ 
and $\vN$ be a unit normal vector field of the hypersurface $M$.  Denote by $\widetilde R$ the curvature tensor of $\widetilde g$.  
Also, denote by $A$ and $\widetilde R(\vN)$ the shape tensor and the normal Jacobi field of the submanifold $M$ for $\vN$, where $\widetilde R(\vN)$ is 
the $(1,1)$-tensor field on $M$ defined by 
$$\widetilde R(\vN)_x:=df_x^{-1}\circ\widetilde R_{f(x)}(\cdot,\vN_x)\vN_x\circ df_x\quad\,\,(x\in M).$$
If, for each $x\in M$, $A_x$ and $\widetilde R(\vN)_x$ commutes (i.e., $[A_x,\widetilde R(\vN)_x]=0$) for each $x\in M$, then the hypersurface $M$ 
is said to be {\it curavture-adapted}.  This notion was introduced by Berndt-Vanhecke (\cite{BV}), where we note that they (\cite{BV}) defined 
the notion of a curvature-adapted submanifold (of general codimension) in more general.  

Let $M$ be a connected oriented hypersurface in a K$\ddot{\rm a}$hler manifold $(\widetilde M,\widetilde g,\widetilde J)$ immersed by 
an immersion $f:M\hookrightarrow\widetilde M$ and $\vN$ be a unit normal vector field of the hypersurface $M$.  
Set $\xi:=\widetilde J(\vN)$.  This tangent vector field $\xi$ is called the {\it structure vector field} of the hypersurface $M$.  
If the integral curves of $\xi$ are geodesics in $M$ (or equivalently, $\xi_x$ is the eignevector of $A_x$ for each $x\in M$), 
then it is called a {\it Hopf hypersurface}.  
For any (real) hypersurface $M$ (immersed by an immersion $f$) in the $n$-dimensional complex projective space $\mathbb{CP}^n(c)(=SU(n+1)/S(U(1)\times U(n)))$ 
of constant holomofrphic sectional curvature $c(>0)$ and the $n$-dimensional complex hyperbolic space $\mathbb{CH}^n(-c)(=SU(1,n)/S(U(1)\times U(n)))$ 
of constant holomorphic sectional curvature $-c(<0)$, the eigenspace decomposition of the normal Jacobi operator $\widetilde R(\vN)_x$ is given by 
$T_xM={\rm Span}\{\xi_x\}\oplus{\rm Span}\{\xi_x\}^{\perp}$.  
From this fact, we can derive directly that $M$ is a Hopf hypersurface if and only if it is curvature-adapted.  
It is known that any Hopf hypersurface in $\mathbb{CP}^n(c)$ is a tube (of constant radius) over a complex submanifold in $\mathbb{CP}^n(c)$ (see \cite{CR}).  
Also, it is known that any Hopf hypersurface in $\mathbb{CH}^n(-c)$ the absolute value of whose principal curvature corresponding to $\xi$ is greater than 
$\sqrt{-c}$ is a tube (of constant radius) over a complex submanifold in $\mathbb{CH}^n(-c)$ (see \cite{Mo}), 
where we note that the complex submanifold is the focal submanifold for the principal curvature corresponding to $\xi$.  
In the case where a Hopf hypersurface $M$ in $\mathbb{CH}^n(-c)$ the absolute value of whose principal curvature corresponding to $\xi$ 
is smaller than $\sqrt{-c}$, focal submanifolds of $M$ vanish beyond the ideal boundary of $\mathbb{CH}^n(-c)$.  
Hence it does not occur as a tube (of constant radius) over a complex submanifold in $\mathbb{CH}^n(-c)$.  
So we need to consider the extrinsic complexification $M^{\mathbb C}$ of $M$.  This extrinsic complexification $M^{\mathbb C}$ is defined as a complex submanifold in the complexification 
$(\mathbb{CH}^n(-c))^{\mathbb C}(=SL(n+1,\mathbb C)/(SL(1,\mathbb C)\times SL(n,\mathbb C)))$ of $\mathbb{CH}^n(-c)$.  
Note that $(\mathbb{CH}^n(-c))^{\mathbb C}$ is a semi-simple pseudo-Riemannian symmetric space with a neutral metric called 
an {\it anti-K$\ddot{\rm a}$hler symmetric space}.  
This extrinsic complexification $M^{\mathbb C}$ admits the focal submanifold of $M^{\mathbb C}$ (in $(\mathbb {CH}^n(-c))^{\mathbb C}$) 
for the principal curvature corresponding to $\xi$ and hence $M^{\mathbb C}$ occurs as a tube (of constant complex radius) 
over the focal submanifold.  See \cite{Koi4} (\cite{Koi3} also) about the notions of the extrinsic complexification 
$M^{\mathbb C}(\subset\widetilde M^{\mathbb C}$) of a real analytic submanifold $M$ in a pseudo-Riemannian manifold $\widetilde M$.  
Also, see \cite{Koi3}, \cite{Koi4} and \cite{Koi5} about the research by using the notion of this extrinsic complexification.  

Let $M$ be a connected oriented hypersurface in the Riemannian product manifold $(\widetilde M_1\times\widetilde M_2,\widetilde g_1\times\widetilde g_2)$ of 
Riemannian manifolds $(\widetilde M_i,\widetilde g_i)$ ($i=1,2$) and $\vN$ be the unit normal vector field of $M$.  
For the simplicity, we abbreviate $(\widetilde M_1\times\widetilde M_2,\widetilde g_1\times\widetilde g_2)$ as $(\widetilde M,\widetilde g)$.  
Also, let $P$ denotes the product structure of the product manifold $\widetilde M$, which is defined by 
$$\begin{array}{c}
P_{(p_1,p_2)}(\vv_1,\vv_2):=(\vv_1,-\vv_2)\\
((p_1,p_2)\in \widetilde M_1\times\widetilde M_2,\,\,(\vv_1,\vv_2)\in T_{(p_1,p_2)}(\widetilde M_1\times\widetilde M_2)(=T_{p_1}\widetilde M_1\oplus T_{p_2}\widetilde M_2)).
\end{array}$$
Set $C:=\widetilde g(P\vN,\vN)$, which is called the {\it product angle function} of this hypersurface $M$.  
Denote by $\mathbb S^2(1)$ (resp. $\mathbb H^2(-1)$) the sphere of constant sectional curvature $1$ (resp. the hyperbolic space of constant curvature $-1$).  
We consider the case where the ambient space $(\widetilde M,\widetilde g)$ is $\mathbb S^2(1)\times\mathbb S^2(1)$ or $\mathbb H^2(-1)\times\mathbb H^2(-1)$.  
Urbano (\cite{U}) classified locally the following hypersurfaces in $\mathbb S^2(1)\times\mathbb S^2(1)$;

\vspace{0.15truecm}

\noindent
$\bullet\,\,$ 
Isoparametric hypersurfaces;

\noindent
$\bullet\,\,$ 
Hypersurfaces with at most two distinct constant principal curvatures;

\noindent
$\bullet\,\,$ 
Hypersurfaces with three distinct constant principal curvatures such that $\vV:=(P\vN)_T$ is a principal vector field, where $(P\vN)_T$ 
denotes the tangent component of $P\vN$.  

\vspace{0.15truecm}

\noindent
On the other hand, Gao-Ma-Yao (\cite{GMY}) classified locally the above hypersurfaces in $\mathbb H^2(-1)\times\mathbb H^2(-1)$.  
According to these local classifications, all the above hypersurfaces are curvature-adapted.  
Thus good hypersurfaces in $\mathbb S^2(1)\times\mathbb S^2(1)$ and $\mathbb H^2(-1)\times\mathbb H^2(-1)$ are curvature-adapted.  
Recently, Hu-Lu-Yao-Zhang (\cite{HLYZ}) classified curvature-adapted hypersurfaces in $\mathbb S^2(1)\times\mathbb S^2(1)$ and 
$\mathbb H^2(-1)\times\mathbb H^2(-1)$.  

Let $G/K$ be a (Riemannian) symmetric space and $H$ be a symmetric subgroup of $G$ (i.e., there exists an involution $\sigma$ 
of $G$ with $({\rm Fix}\,\sigma)_0\subset H\subset{\rm Fix}\,\sigma$), where ${\rm Fix}\,\sigma$ is the fixed point set of $\sigma$ 
and $({\rm Fix}\,\sigma)_0$ is the identiy component of ${\rm Fix}\,\sigma$.  
Then the natural action $H\curvearrowright G/K$ is called a {\it Hermann action} in the case where $G/K$ is of compact type 
and it is called a {\it Hermann type action} in the case where $G/K$ is of non-compact type.  
It is known that the principal orbits of this action are curvature-adapted submanifolds.  This fact was proved in \cite{GT} 
in the case where $G/K$ is of compact type and that it was proved in \cite{Koi1} in the case where $G/K$ is of non-compact type 
(in the method of the proof different from that of \cite{GT}).  
See \cite{Koi2} about the research of Hermann type actions on a pseudo-Riemannian symmetric space in more general.  

In this paper, we give a method constructing complete curvature-adapted hypersurfaces in the product $G_1/K_1\times G_2/K_2$ 
from a pair of complete curvature-adapted hypersurfaces in $G_1/K_1$ and $G_2/K_2$.  Let $M_1$ and $M_2$ be curvature-adapted 
hypersurfaces in symmetric spaces $G_1/K_1$ and $G_2/K_2$ immersed by immersions $f_i:M_i\hookrightarrow G_i/K_i$, respectively.  
Denote by $\vN_i$ ($i=1,2$) a unit normal vector field of the hypersurface $M_i$.  Let $S^1$ be the unit circle, that is, 
$S^1:=\{(x,y)\,|\,x^2+y^2=1\}$.  Let $\vu=(u_1,u_2):[0,2\pi]\to\mathbb R^2$ be a simple closed regular curve of $C^{\infty}$-class 
surrounding $(0,0)$ such that $\displaystyle{\mathop{\max}_{\theta\in[0,2\pi]}\,\|\vu(\theta)\|}$ 
is sufficiently small, where ``regular'' means that $\vu'(\theta)\not={\bf 0}$ at each $\theta\in[0,2\pi]$ and 
we assume $\vu'(0)=\vu'(2\pi)$.  Also, let $\iota:[0,2\pi)\to S^1$ be the bijection defined by 
$\iota(\theta):=(\cos\theta,\sin\theta)$ ($\theta\in[0,2\pi)$).  
Define a map $f:M_1\times M_2\times S^1\to G_1/K_1\times G_2/K_2$ by 
$$\begin{array}{r}
\displaystyle{f(p_1,p_2,\iota(\theta)):=(\exp_1(u_1(\theta)(\vN_1)_{p_1}),\exp_2(u_2(\theta)(\vN_2)_{p_2}))}\\
\displaystyle{(p_1,p_2,\iota(\theta))\in M_1\times M_2\times S^1),}
\end{array}$$
where $\exp_i$ ($i=1,2$) denote the normal exponential maps of the hypersurfaces $M_i$ ($i=1,2$).  
The main theorem in this paper is as follows.  

\vspace{0.25truecm}

\noindent
{\bf Theorem A.}\ {\sl If $\displaystyle{\mathop{\max}_{\theta\in[0,2\pi]}\,\|\vu(\theta)\|}$ is sufficiently small, then the above map $f$ is 
an immersion and $M_1\times M_2\times S^1$ is a curvature-adapted hypersurface in $G_1/K_1\times G_2/K_2$ immersed by $f$.  
The eigenvalues of the shape operator $A$ and the normal Jacobi operator $\widetilde R(\overline{\vN})$ of this curvature-adapted hypersurface 
are as in Table 1, where $\overline{\vN}$ is the unit normal vector field of $f$.  Also, the product angle function $C$ of this curvature-adapted 
hypersurface is given by $\displaystyle{C(p_1,p_2,\iota(\theta))=-\frac{u_1'(\theta)^2-u_2'(\theta)^2}{u_1'(\theta)^2+u_2'(\theta)^2}}$.}

\vspace{0.25truecm}

\newpage

$$\begin{tabular}{|c|c|}
\hline
{\scriptsize Common eigenspaces} & {\scriptsize Eigenvalues of $A_{(p_1,p_2,\iota(\theta))}$}\\
\hline
{\scriptsize $E^1_{jk,p_1}$} & {\scriptsize $-\frac{u'_2(\theta)}{\sqrt{u'_1(\theta)^2+u'_2(\theta)^2}}\cdot
\frac{\sqrt{\mu_{1k,p_1}}\sin(u_1(\theta)\sqrt{\mu_{1k,p_1}})+\lambda_{1j,p_1}\cos(u_1(\theta)\sqrt{\mu_{1k,p_1}})}
{\cos(u_1(\theta)\sqrt{\mu_{1k,p_1}})-\frac{\lambda_{1j,p_1}\sin(u_1(\theta)\sqrt{\mu_{1k,p_1}})}{\sqrt{\mu_{1k,p_1}}}}$}\\
\hline
{\scriptsize $E^2_{jk,p_2}$} & {\scriptsize $\frac{u'_1(\theta)}{\sqrt{u'_1(\theta)^2+u'_2(\theta)^2}}\cdot
\frac{\sqrt{\mu_{2k,p_2}}\sin(u_2(\theta)\sqrt{\mu_{2k,p_2}})+\lambda_{2j,p_2}\cos(u_2(\theta)\sqrt{\mu_{2k,p_2}})}
{\cos(u_2(\theta)\sqrt{\mu_{2k,p_2}})-\frac{\lambda_{2j,p_2}\sin(u_2(\theta)\sqrt{\mu_{2k,p_2}})}{\sqrt{\mu_{2k,p_2}}}}$}\\
\hline
{\scriptsize ${\rm Span}\{\frac{\partial}{\partial\theta}\}$} & {\scriptsize $\frac{u'_1(\theta)u''_2(\theta)-u''_1(\theta)u'_2(\theta)}
{(u'_1(\theta)^2+u'_2(\theta)^2)^{\frac{3}{2}}}$}\\
\hline
\end{tabular}$$

$$\begin{tabular}{|c|c|}
\hline
{\scriptsize Common eigenspaces} & {\scriptsize Eigenvalues of $\widetilde R(\overline{\vN})_{(p_1,p_2,\iota(\theta))}$}\\
\hline
{\scriptsize $E^1_{jk,p_1}$} & {\scriptsize $\frac{u'_2(\theta)^2\mu_{1k,p_1}}{u'_1(\theta)^2+u'_2(\theta)^2}$}\\
\hline
{\scriptsize $E^2_{jk,p_2}$} & {\scriptsize $\frac{u'_1(\theta)^2\mu_{2k,p_2}}{u'_1(\theta)^2+u'_2(\theta)^2}$}\\
\hline
{\scriptsize ${\rm Span}\{\frac{\partial}{\partial\theta}\}$} & {\scriptsize $0$}\\
\hline
\end{tabular}$$

\vspace{0truecm}

{\scriptsize 
$$\left(\begin{array}{l}
T_{p_i}M_i=\oplus_{(j,k)\in I_{A_iR_i,p_i}}E^i_{jk,p_i}\,:\,{\rm the}\,\,{\rm common}\,\,{\rm eigenspace}\,\,{\rm decomposition}\,\,{\rm of}\,\,{\rm the}\,\,{\rm shape}\,\,{\rm operator}
\,\,(A_i)_{p_i}\\
\hspace{4.45truecm}\,\,{\rm and}\,\,{\rm the}\,\,{\rm normal}\,\,{\rm Jacobi}\,\,{\rm operator}\,\,\widetilde R_i(\vN_i)_{p_i}\,\,{\rm of}\,\,M_i\\
\lambda_{ij,p_i}\,:\,{\rm the}\,\,{\rm eigenvalue}\,\,{\rm of}\,\,(A_i)_{p_i}\,\,{\rm corresponding}\,\,{\rm to}\,\,\oplus_kE^i_{jk,p_i}\\
\mu_{ik,p_i}\,:\,{\rm the}\,\,{\rm eigenvalue}\,\,{\rm of}\,\,R_i(\vN_i)_{p_i}\,\,{\rm corresponding}\,\,{\rm to}\,\,\oplus_jE^i_{jk,p_i}
\end{array}\right)$$
}

\vspace{0.25truecm}

\centerline{{\bf Table 1:$\,\,$Eigenvalues of $A_{(p_1,p_2,\iota(\theta))}$ and $\widetilde R(\overline{\vN})_{(p_1,p_2,\iota(\theta))}$}}

\vspace{0.5truecm}

\noindent
{\it Remark 1.1.}\ \ Denote by $T_{(p_1,p_2)}^{\perp}(M_1\times M_2)$ the normal space of the the submanifold $(f_1\times f_2)(M_1\times M_2)$ 
of codimesnsion two at $(p_1,p_2)$ and set $\Sigma_{(p_1,p_2)}:=\exp(T_{(p_1,p_2)}^{\perp}(M_1\times M_2))$, where $\exp$ denotes 
the exponential map of $G_1/K_1\times G_2/K_2$.  It is easy to show that the induced metric on $\Sigma_{(p_1,p_2)}$ is flat and hence 
$\Sigma_{(p_1,p_2)}$ is totally geodesic in $G_1/K_1\times G_2/K_2$.  From this fact, we see that $(f_1\times f_2)(M_1\times M_2)$ is 
a submanifold with flat section in $G_1/K_1\times G_2/K_2$.  
Since $f(\{(p_1,p_2)\}\times S^1)$ is a simple closed regular curve surrounding $(f_1\times f_2)(p_1,p_2)$ in the flat section $\Sigma_{(p_1,p_2)}$, 
$f(M_1\times M_2\times S^1)$ is a tubular-like hypersurface over $(f_1\times f_2)(M_1\times M_2)$.  
In particular, in the case where $\vu$ is given by $\vu(\theta):=(r\cos\theta,r\sin\theta)$ ($\theta\in[0,2\pi]$) for some positive (small) constant $r$, 
$f(M_1\times M_2\times S^1)$ is the tube of radius $r$ over $(f_1\times f_2)(M_1\times M_2)$.  

\section{Basic notions and facts} 
\subsection{Curvature-adapted hypersurfaces} 
Let $M$ be a connected oriented hypersurface in a Riemannian manifold $(\widetilde M,\widetilde g)$ immersed by an immersion $f:M\hookrightarrow\widetilde M$ 
and $\vN$ be a unit normal vector field of this hypersurface $M$.  
Denote by $\widetilde{\nabla}$ and $\widetilde R$ the Riemannian connection and the curvature tensor of $\widetilde g$, respectively.  
Also, denote by $A$ and $\widetilde R(\vN)$ the shape tensor and the normal Jacobi operator of the submanifold $M$ for $\vN$, where 
$\widetilde R(\vN)$ is the $(1,1)$-tensor field on $M$ defined by 
$$\widetilde R(\vN)_p:=df_p^{-1}\circ\widetilde R_{f(p)}(\cdot,\vN_p)\vN_p\circ df_p\quad\,\,(p\in M).$$
If, for each $p\in M$, $A_p$ and $\widetilde R(\vN)_p$ commutes (i.e., $[A_p,\widetilde R(\vN)_p]=0$) for each $p\in M$, then 
the hypersurface $M$ is said to be {\it curavture-adapted}.  
This notion was introduced by Berndt-Vanhecke (\cite{BV}).  

Assume that $M$ is curvature-adapted.  Denote by $g$ the induced metric $f^{\ast}\widetilde g$ on $M$.  
Since $A_p$ and $\widetilde R(\vN)_p$ are symmetric $(1,1)$-tensor on $(T_pM,g_p)$, they are diagonalizeble with respect to an orthonormal basis 
and we have the orthogonal decompositions 
$\displaystyle{T_pM=\mathop{\oplus}_{i\in I_{A,p}}E_{i,p}^A}$ cosisting of eigenspaces of $A_p$ and the orthogonal decomposition 
$\displaystyle{T_pM=}$\newline
$\displaystyle{\mathop{\oplus}_{i\in I_{R,p}}E_{i,p}^R}$ consisting of the eigenspace of $\widetilde R(\vN)_p$.  
%Assume that $I_{A,p}$ and $I_{R,p}$ are independent of the choice of $p\in M$.  For the simplicity, we denote $I_{A,p}$ (resp. $I_{R,p}$) by $I_{A,p}$ (resp. $I_R$).  
%Then we can show that the correspondences $x\mapsto(E_i^A)_p$ ($p\in M$) and $x\mapsto(E_i^R)_p$ ($p\in M$) give $C^{\infty}$-distributions on $M$.  
%Denote by $E_i^A$ and $E_i^R$ these distributions.  
Define two families $\{\lambda_{i,p}\}_{i\in I_{A,p}}$ and $\{\mu_{i,p}\}_{i\in I_{R,p}}$ of real numbers by 
$$A_p|_{E_{i,p}^A}=\lambda_{i,p}\,{\rm id}\quad\,\,{\rm and}\quad\,\,\widetilde R(\vN)_p|_{E_{i,p}^R}=\mu_{i,p}\,{\rm id}.$$
%These functions $\lambda_{i,p}$ ($i\in I_{A,p}$) and $\mu_{j,p}$ ($j\in I_{R,p}$) are called the {\it eigenvalues} of $A$ and $\widetilde R(\vN)$, respectively.  
%Also, the distributions $E_i^A$ ($i\in I_A$) and $E_j^R$ ($j\in I_R$) are called the {\it eigendistributions} of $A$ and $\widetilde R(\vN)$, respectively.  
Define linear subspaces $E_{ij,p}$ ($(i,j)\in I_{A,p}\times I_{R,p}$) on $M$ by 
$E_{ij,p}:=E_{i,p}^A\cap E_{j,p}^R$.  Set 
$$I_{AR,p}:=\{(i,j)\in I_{A,p}\times I_{R,p}\,|\,{\rm dim}\,E_{ij,p}\not=0\}.$$
Since $A_p$ and $\widetilde R(\vN)_p$ commute, they have simultaneously daigonalizable, that is, 
the following orthogonal decomposition holds:
$$T_pM=\mathop{\oplus}_{(i,j)\in I_{AR,p}}E_{ij,p}.$$
This orthogonal decomposition is called the {\it common eigenspace decomposition} of $A_p$ and $\widetilde R(\vN)_p$.  

\subsection{Strongly $M$-Jacobi fields} 
In this section, we recall the explicit description of a strongly $M$-Jacobi field for a connected oriented hypersurface $M$ in a symmetric space $G/K$.  
Let $M$ be an immersed submanifold in $G/K$ immersed by an immersion $f$.  Fix a unit normal vector field $\vN$ of $M$.  
Denote by $\nabla,\,\widetilde{\nabla}$ and $A$ the Reimannian connection of the induced metric on $M$, that of the metric of $G/K$ 
and the shape operator of $M$ for $\vN$, respectively.  Also, let $\widetilde R(\vN)$ be the normal Jacobi operator of $M$ for $\vN$.  
Take $p\in M$.  Let $\gamma_p:[0,\infty)\to G/K$ be the geodesic of the direction $\vN_p$ (i.e., $\gamma_p'(0)=\vN_p$) and 
$\vY$ be a Jacobi field along $\gamma_p$, that is, a vector field along $\gamma_p$ satisfying the Jacobi equation:
$$\vY''(s)+\widetilde R_{\gamma_p(s)}(\vY'(s),\gamma_p'(s))\gamma_p'(s)={\bf 0},$$
where $\vY'(s)$ (resp. $\vY''(s)$) denotes $\widetilde{\nabla}^{\gamma}_{\frac{d}{ds}}\vY$ 
(resp. $\widetilde{\nabla}^{\gamma}_{\frac{d}{ds}}(\widetilde{\nabla}^{\gamma}_{\frac{d}{ds}}\vY)$).  
In particular, if $\vY(0)$ and $\vY'(0)$ belong to $df_p(T_pM)$, then $\vY$ is called a {\it strongly $M$-Jacobi field}.  
Then we note that $\vY'(0)=-df_p(A_p(\vY(0)))$ holds.  Since $\widetilde{\nabla}\widetilde R={\bf 0}$, the strongly $M$-Jacobi field $Y$ 
is described as 
$$\vY(s)=P_{\gamma_p\vert_{[0,s]}}\left(df_p\left(\left(\cos\left(s\sqrt{\widetilde R(\vN)_p}\right)
-\frac{\sin\left(s\sqrt{\widetilde R(\vN)_p}\right)}{\sqrt{\widetilde R(\vN)_p}}\circ A_p\right)(\vY(0))\right)\right),\leqno{(2.1)}$$
where $P_{\gamma_p\vert_{[0,s]}}$ is the parallel translation along $\gamma_p\vert_{[0,s]}$, and $\displaystyle{\cos\left(s\sqrt{\widetilde R(\vN)_p}\right)}$ \newline
and $\displaystyle{\sin\left(s\sqrt{\widetilde R(\vN)_p}\right)/\sqrt{\widetilde R(\vN)_p}}$ are defined by 
$$\begin{array}{c}
\displaystyle{\cos\left(s\sqrt{\widetilde R(\vN)_p}\right):=\sum_{j=0}^{\infty}\frac{(-1)^js^{2j}}{(2j)!}\cdot{\widetilde R(\vN)_p}^j,}\\
\displaystyle{\frac{\sin\left(s\sqrt{\widetilde R(\vN)_p}\right)}{\sqrt{\widetilde R(\vN)_p}}:=\sum_{j=0}^{\infty}\frac{(-1)^js^{2j+1}}{(2j+1)!}
\cdot{\widetilde R(\vN)_p}^j,}
\end{array}$$
respectively.  

\section{Proof of Theorem A} 
In this section, we shall prove Theorem A.  

\vspace{0.25truecm}

\noindent
{\it Proof of Theorem A.}\ \ We shall use the notations stated in Introduction.  Denote by $\widetilde{\nabla}^i$ and $\widetilde R_i$ the Riemannian connection 
and the curvature tensor of $G_i/K_i$ ($i=1,2$), and denote by $\widetilde{\nabla}$ and $\widetilde R$ those of $G_1/K_1\times G_2/K_2$.  
Also, denote by $A_i$ ($i=1,2$) the shape opertors of the hypersurface $M_i$ in $G_i/K_i$ immersed by $f_i$.  
Fix $(p_1,p_2,\iota(\theta))\in M_1\times M_2\times S^1$.  
Let $\{E_{j,p_i}^{A_i}\}_{j\in I_{A_i,p_i}}$ and $\{E_{j,p_i}^{R_i}\}_{j\in I_{R_i,p_i}}$ be the families of eigenspaces of $(A_i)_{p_i}$ 
and $\widetilde R_i(\vN_i)_{p_i}$.  
Also, let $\lambda_{ij,p_i}$ be the eignvalue of $(A_i)_{p_i}$ corresponding to $E_{j,p_i}^{A_i}$ and 
$\mu_{ij,p_i}$ be the eignvalue of $\widetilde R_i(\vN_i)_{p_i}$ corresponding to $E_{j,p_i}^{R_i}$.  
%Let $\displaystyle{T_xM_i=\mathop{\oplus}_{j\in I_{A_i}^x}(E_j^{A_i})_x}$ be the orthogonal decomposition cosisting of eigenspaces of $(A_i)_x$ 
%and the orthogonal decomposition $\displaystyle{T_xM_i=\mathop{\oplus}_{j\in I_{R_i}^x}(E_j^{R_i})_x}$ consisting of the eigenspace of $\widetilde R(\vN_i)_x$.  
%For the simplicity, we consider the case where $I_{A_i}^x$ and $I_{R_i}^x$ are independent of the choice of $x\in M_i$.  
%For the simplicity, we denote $I_{A_i}^x$ (resp. $I_{R_i}^x$) by $I_{A_i}$ (resp. $I_{R_i}$).  
%As stated in Subsection 2.2, the $C^{\infty}$-distributions $E_j^{A_i}$ and $E_j^{R_i}$ on $M_i$  are defined.  
%Also, two families $\{\lambda_{ij}\}_{j\in I_{A_i}}$ and $\{\mu_{ij}\}_{j\in I_{R_i}}$ of $C^{\infty}$-functions on $M_i$ by 
%$$(A_i)_x|_{(E_j^{A_i})_x}=\lambda_{ij}(x)\,{\rm id}\quad\,\,{\rm and}\quad\,\,\widetilde R_i(\vN_i)_x|_{(E_j^{R_i})_x}=\mu_{ij}(x)\,{\rm id}.$$
%These functions $\lambda_{ij}$ ($j\in I_{A_i}$) and $\mu_{ij}$ ($j\in I_{R_i}$) are called the {\it eigenfunctions} of $A_i$ and $\widetilde R_i(\vN_i)$, respectively.  
%Also, the distributions $E_j^{A_i}$ ($j\in I_{A_i}$) and $E_j^{R_i}$ ($j\in I_{R_i}$) are called the {\it eigendistributions} of $A_i$ and $\widetilde R_i(\vN_i)$, respectively.  
Define linear subspaces $E_{jk,p_i}^i$ ($(j,k)\in I_{A_i,p_i}\times I_{R_i,p_i}$) on $M$ by 
$E_{jk,p_i}^i:=E_{j,p_i}^{A_i}\cap E_{k,p_i}^{R_i}$.  Set 
$$I_{A_iR_i,p_i}:=\{(j,k)\in I_{A_i,p_i}\times I_{R_i,p_i}\,|\,{\rm dim}\,E_{jk,p_i}^i\not=0\}.$$
Since $(A_i)_{p_i}$ and $\widetilde R_i(\vN_i)_{p_i}$ commute by the assumption, they have simultaneously daigonalizable, that is, 
the following orthogonal decomposition holds:
$$T_{p_i}M_i=\mathop{\oplus}_{(j,k)\in I_{A_iR_i,p_i}}E_{jk,p_i}^i.\leqno{(3.1)}$$
Hence we have 
$$T_{(p_1,p_2)}(M_1\times M_2)=T_{p_1}M_1\oplus T_{p_2}M_2=\left(\mathop{\oplus}_{(j,k)\in I_{A_1R_1,p_1}}E_{jk,p_1}^1\right)\oplus
\left(\mathop{\oplus}_{(j,k)\in I_{A_2R_2,p_2}}E_{jk,p_2}^2\right).\leqno{(3.2)}$$

Take a unit normal vector field $\overline{\vN}$ of the hypersurface $M_1\times M_2\times S^1$ in $G_1/K_1\times G_2/K_2$ immersed by $f$.  
%%, where $\overline{\vN}_i$ denotes the unit normal vector fireld of the hypersurface $M_i$ in $G_I/K_i$ immersed by $f_i$ and 
%%$\rho_i$ ($i=1,2$) are positive functions over $M_1\times M_2\times S^1$ with $\rho_1^2+\rho_2^2=1$.  
%%For the convenience, set $\overline{\vN}_i:=\rho_i\vN_i$ ($i=1,2$).  
Note that 
\begin{align*}
&\overline{\vN}_{(p_1,p_2,\iota(\theta))}=(\rho_1(p_1,p_2,\iota(\theta))\,\overline{\vN_1}_{p_1},\,\rho_2(p_1,p_2,\iota(\theta))\,\overline{\vN_2}_{p_2})\\
\in&T_{(p_1,p_2,\iota(\theta))}^{\perp}(M_1\times M_2\times S^1)\\
\subset&T_{f(p_1,p_2,\iota(\theta))}(G_1/K_1\times G_2/K_2)\\
=&\mathop{\oplus}_{i=1}^2T_{\exp_i(u_i(\theta)(\vN_i)_{p_i})}(G_i/K_i).
\end{align*}
Denote by $A$ the shape operator of the hypersurface $M_1\times M_2\times S^1$ in $G_1/K_1\times G_2/K_2$ (immersed by $f$) for $\overline{\vN}$.  
For a positive constant $r$, define a map $(f_i)_r:M_i\to G_i/K_i$ by $(f_i)_r(p):=\exp_i(r(\vN_i)_p)$ ($p\in M_i$) and set $(M_i)_r:=(f_i)_r(M_i)$.  

Take $\vv^i_{jk}\in E^i_{jk,p_i}$.  Let $c^i_{jk}:(-\varepsilon,\varepsilon)\to G_i/K_i$ be a curve with $(c^i_{jk})'(0)=\vv^i_{jk}$.  Then we have 
$$\begin{array}{l}
\hspace{0.5truecm}\displaystyle{df_{(p_1,p_2,\iota(\theta))}(\vv^1_{jk})=(f\circ c^1_{jk})'(0)}\\
\displaystyle{=\left.\frac{d}{dt}\right|_{t=0}\left(\exp_1(u_1(\theta)(\vN_1)_{c^1_{jk}(t)}),\,\exp_2(u_2(\theta)(\vN_2)_{p_2})\right)}\\
\displaystyle{=((d(f_1)_{u_1(\theta)})_{p_1}(\vv^1_{jk}),{\bf 0})\,\,(\in d(f_1)_{u_1(\theta)})_{p_1}(E^1_{jk,p_1})).}
\end{array}\leqno{(3.3)}$$
Similarly, we have 
$$\begin{array}{l}
\hspace{0.5truecm}\displaystyle{df_{(p_1,p_2,\iota(\theta))}(\vv^2_{jk})=(f\circ c^2_{jk})'(0)}\\
\displaystyle{=\left.\frac{d}{dt}\right|_{t=0}\left(\exp_1(u_1(\theta)(\vN_1)_{p_1}),\,\exp_2(u_2(\theta)(\vN_2)_{c^2_{jk}(t)})\right)}\\
\displaystyle{=({\bf 0},(d(f_2)_{u_2(\theta)})_{p_2}(\vv^2_{jk}))\,\,(\in d((f_2)_{u_2(\theta)})_{p_2}(E^2_{jk,p_2})).}
\end{array}\leqno{(3.4)}$$
Let $\gamma_{p_i}:[0,\infty)\to G_i/K_i$ be the geodesic in $G_i/K_i$ of direction $(\vN_i)_{p_i}$.  
Also, we have 
$$\begin{array}{l}
\hspace{0.5truecm}\displaystyle{df_{(p_1,p_2,\iota(\theta))}\left(\left(\frac{\partial}{\partial\theta}\right)_{(p_1,p_2,\iota(\theta))}\right)}\\
\displaystyle{=\frac{d}{d\theta}\left(\exp_1(u_1(\theta)(\vN_1)_{p_1}),\,\exp_2(u_2(\theta)(\vN_2)_{p_2})\right)}\\
\displaystyle{=\left(u_1'(\theta)P_{\gamma_{p_1}|_{[0,u_1(\theta)]}}((\vN_1)_{p_1}),\,u_2'(\theta)P_{\gamma_{p_2}|_{[0,u_2(\theta)]}}((\vN_2)_{p_2})\right)}\\
\displaystyle{(\in {\rm Span}\{P_{\gamma_{p_1}|_{[0,u_1(\theta)]}}((\vN_1)_{p_1})\}\oplus{\rm Span}\{P_{\gamma_{p_2}|_{[0,u_2(\theta)]}}((\vN_2)_{p_2})\}).}
\end{array}\leqno{(3.5)}$$
For the simplicity, denote by $c_i$ the restricted curve $c^i_{jk}|_{(-\varepsilon,\varepsilon)}$, wehre $\varepsilon$ is a sufficiently small positive number.  
Define a map $\delta_i:(-\varepsilon,\varepsilon)\times[0,\infty)\to G_i/K_i$ by 
$$\delta_i(t,s):=\exp_i(s(\vN_i)_{c_i(t)})\quad((t,s)\in(-\varepsilon,\varepsilon)\times[0,\infty))$$
and a vector field $\vY_i$ along $\gamma_{p_i}$ by $\vY_i(s):=(d\delta_i)_{(0,s)}\left(\left(\frac{\partial}{\partial t}\right)_{(0,s)}\right)$.  
Since $\delta_i$ is a geodesi variation consisting of unit speed normal geodesics of $M_i$, its variational vector field $\vY_i$ is a strongly $M_i$-Jacobi field.  
Hence, according to $(2.1)$, we have 
$$\begin{array}{l}
\hspace{0.5truecm}\displaystyle{\vY_i(s)}\\
\displaystyle{=P_{\gamma_{p_i}|_{[0,s]}}\left((df_i)_{p_i}\left(\left(\cos\left(s\sqrt{\widetilde R_i(\vN_i)_{p_i})}\right)
-\frac{\sin\left(s\sqrt{\widetilde R_i(\vN_i)_{p_i})}\right)}{\sqrt{\widetilde R_i(\vN_i)_{p_i})}}\circ(A_i)_{p_i}\right)(\vv^i_{jk})\right)\right)}\\
\displaystyle{=\left(\cos(s\sqrt{\mu_{ik,p_i}})-\frac{\lambda_{ij,p_i}\sin(s\sqrt{\mu_{ik,p_i}})}{\sqrt{\mu_{ik,p_i}}}\right)P_{\gamma_{p_i}|_{[0,s]}}
((df_i)_{p_i}(\vv^i_{jk}))}\\
\displaystyle{\in P_{\gamma_{p_i}|_{[0,s]}}((df_i)_{p_i}(E^i_{jk,p_1})).}
\end{array}\leqno{(3.6)}$$
Here, in the case of $\mu_{ik,p_i}<0$, $\cos(s\sqrt{\mu_{ik,p_i}})$ and $\frac{\sin(s\sqrt{\mu_{ik,p_i}})}{\sqrt{\mu_{ik,p_i}}}$ imply 
$\cosh(s\sqrt{-\mu_{ik,p_i}})$ and $\frac{\sinh(s\sqrt{-\mu_{ik,p_i}})}{\sqrt{-\mu_{ik,p_i}}}$, respectively.  
On the other hand, from the definition of $\delta_i$,  we have $\vY_i(u_i(\theta))=d((f_i)_{u_i(\theta)})_{p_i}(\vv^i_{jk})$.  
Hence we obtain 
$$\begin{array}{l}
\hspace{0.5truecm}\displaystyle{d((f_i)_{u_i(\theta)})_{p_i}(\vv^i_{jk})}\\
\displaystyle{=\left(\cos(u_i(\theta)\sqrt{\mu_{ik,p_i}})-\frac{\lambda_{ij}(p_i)\sin(u_i(\theta)\sqrt{\mu_{ik,p_i}})}{\sqrt{\mu_{ik,p_i}}}\right)P_{\gamma_{p_i}|_{[0,u_i(\theta)]}}
((df_i)_{p_i}(\vv^i_{jk}))}\\
\displaystyle{\in P_{\gamma_{p_i}|_{[0,u_i(\theta)]}}(d(f_i)_{p_i}(E^i_{jk,p_i})).}
\end{array}\leqno{(3.7)}$$
From $(3.3),\,(3.4)$ and $(3.7)$, we can derive 
$$\begin{array}{l}
\hspace{0.5truecm}\displaystyle{df_{(p_1,p_2,\iota(\theta))}(\vv^1_{jk})}\\
\displaystyle{\left(\,
\left(\cos(u_1(\theta)\sqrt{\mu_{1k}(p_1)})-\frac{\lambda_{1j,p_1}\sin(u_1(\theta)\sqrt{\mu_{1k,p_1}})}{\sqrt{\mu_{1k,p_1}}}\right)P_{\gamma_{p_1}|_{[0,u_1(\theta)]}}
((df_1)_{p_1}(\vv^1_{jk})),\,{\bf 0}\,\right)}\\
\displaystyle{\in P_{\gamma_{p_1}|_{[0,u_1(\theta)]}}(d(f_1)_{p_1}(E^1_{jk,p_1}))\oplus{\rm Span}\{{\bf 0}\}\,(\subset T_{\gamma_{p_1}(u_1(\theta))}(G_1/K_1)\oplus 
T_{\gamma_{p_2}(u_2(\theta))}(G_2/K_2))}
\end{array}\leqno{(3.8)}$$
and
$$\begin{array}{l}
\hspace{0.5truecm}\displaystyle{df_{(p_1,p_2,\iota(\theta))}(\vv^2_{jk})}\\
\displaystyle{\left(\,{\bf 0},\,\,
\left(\cos(u_2(\theta)\sqrt{\mu_{2k,p_2}})-\frac{\lambda_{2j,p_2}\sin(u_2(\theta)\sqrt{\mu_{2k,p_2}})}{\sqrt{\mu_{2k,p_2}}}\right)P_{\gamma_{p_2}|_{[0,u_2(\theta)]}}
((df_2)_{p_2}(\vv^2_{jk}))\,\right)}\\
\displaystyle{\in{\rm Span}\{{\bf 0}\}\oplus P_{\gamma_{p_2}|_{[0,u_2(\theta)]}}(d(f_2)_{p_2}(E^2_{jk,p_2}))
(\subset T_{\gamma_{p_1}(u_1(\theta))}(G_1/K_1)\oplus T_{\gamma_{p_2}(u_2(\theta))}(G_2/K_2)).}
\end{array}\leqno{(3.9)}$$
Since $\displaystyle{\mathop{\max}_{\theta\in[0,2\pi]}\|\vu(\theta)\|}$ is sufficiently small by the assumption, we obtain 
$$\cos(u_i(\theta)\sqrt{\mu_{ik,p_i}})-\frac{\lambda_{ij,p_i}\sin(u_i(\theta)\sqrt{\mu_{ik,p_i}})}{\sqrt{\mu_{ik,p_i}}}\not=0\quad\,\,(i=1,2).$$
This fact together with $(3.5),\,(3.8)$ and $(3.9)$ implies that $f$ is an immersion.  
Also, from $(3.8)$ and $(3.9)$, we have 
$$\begin{array}{l}
\hspace{0.5truecm}\displaystyle{df_{(p_1,p_2,\iota(\theta))}\left(E^1_{jk,p_1}\oplus{\rm Span}\{{\bf 0}\}\oplus{\rm Span}\{{\bf 0}\}\right)}\\
\displaystyle{=P_{\gamma_{p_1}|_{[0,u_1(\theta)]}}\left((df_1)_{p_1}(E^1_{jk,p_1})\right)\oplus{\rm Span}\{{\bf 0}\}}
\end{array}\leqno{(3.10)}$$
and 
$$\begin{array}{l}
\hspace{0.5truecm}\displaystyle{df_{(p_1,p_2,\iota(\theta))}({\rm Span}\{{\bf 0}\}\oplus E^2_{jk,p_2}\oplus{\rm Span}\{{\bf 0}\})}\\
\displaystyle{={\rm Span}\{{\bf 0}\}\oplus P_{\gamma_{p_2}|_{[0,u_2(\theta)]}}\left((df_2)_{p_2}(E^2_{jk,p_2})\right).}
\end{array}\leqno{(3.11)}$$
Furthermore, from $(3.10)$ and $(3.11)$, we have 
$$\begin{array}{l}
\hspace{0.5truecm}\displaystyle{df_{(p_1,p_2,\iota(\theta))}\left(T_{p_1}M_1\oplus{\rm Span}\{{\bf 0}\}\oplus{\rm Span}\{{\bf 0}\}\right)}\\
\displaystyle{=P_{\gamma_{p_1}|_{[0,u_1(\theta)]}}\left((df_1)_{p_1}(T_{p_1}M_1)\right)\oplus{\rm Span}\{{\bf 0}\}}
\end{array}\leqno{(3.12)}$$
and 
$$\begin{array}{l}
\hspace{0.5truecm}\displaystyle{df_{(p_1,p_2,\iota(\theta))}({\rm Span}\{{\bf 0}\}\oplus T_{p_2}M_2\oplus{\rm Span}\{{\bf 0}\})}\\
\displaystyle{={\rm Span}\{{\bf 0}\}\oplus P_{\gamma_{p_2}|_{[0,u_2(\theta)]}}\left((df_2)_{p_2}(T_{p_2}M_2)\right).}
\end{array}\leqno{(3.13)}$$
According to $(3.12)$ and $(3.13)$, the unit normal vector field $\overline{\vN}$ is described as 
\begin{align*}
\overline{\vN}_{(p_1,p_2,\iota(\theta))}=&\left(\rho_1(p_1,p_2,\iota(\theta))P_{\gamma_{p_1}|_{[0,u_1(\theta)]}}((\vN_1)_{p_1}),\right.\\
&\quad\,\left.\rho_2(p_1,p_2,\iota(\theta))P_{\gamma_{p_2}|_{[0,u_2(\theta)]}}((\vN_2)_{p_2})\right)
\end{align*}
for some functions $\rho_i$ ($i=1,2$) on $M_1\times M_2\times S^1$.  
We shall express these functions in terms of $u_i$ ($i=1,2$).  
By using $(3.5)$, we can derive 
$$\begin{array}{l}
\displaystyle{0=\widetilde g_{f(p_1,p_2,\iota(\theta))}\left(\overline{\vN}_{(p_1,p_2,\iota(\theta))},df_{(p_1,p_2,\iota(\theta))}
\left(\left(\frac{\partial}{\partial\theta}\right)_{(p_1,p_2,\iota(\theta))}\right)\right)}\\
\hspace{0.3truecm}\displaystyle{=\rho_1(p_1,p_2,\iota(\theta))\,u'_1(\theta)+\rho_2(p_1,p_2,\iota(\theta))\,u'_2(\theta).}
\end{array}$$
On the othe hand, we have 
$$1=\widetilde g_{f(p_1,p_2,\iota(\theta))}\left(\overline{\vN}_{(p_1,p_2,\iota(\theta))},\,\overline{\vN}_{(p_1,p_2,\iota(\theta))}\right)
=\rho_1(p_1,p_2,\iota(\theta))^2+\rho_2(p_1,p_2,\iota(\theta))^2.$$
From these relations, we obtain 
$$(\rho_1(p_1,p_2,\iota(\theta)),\rho_2(p_1,p_2,\iota(\theta)))=\frac{1}{\sqrt{u'_1(\theta)^2+u'_2(\theta)^2}}\,(-u'_2(\theta),u'_1(\theta))$$
and hence 
$$\begin{array}{l}
\displaystyle{\overline{\vN}_{(p_1,p_2,\iota(\theta))}=\frac{1}{\sqrt{u'_1(\theta)^2+u'_2(\theta)^2}}
\left(-u'_2(\theta)\,P_{\gamma_{p_1}|_{[0,u_1(\theta)]}}((\vN_1)_{p_1}),\,\right.}\\
\hspace{5.85truecm}\displaystyle{\left.u'_1(\theta)\,P_{\gamma_{p_2}|_{[0,u_2(\theta)]}}((\vN_2)_{p_2})\right).}
\end{array}\leqno{(3.14)}$$
Hence the product angle function $C$ of this hypersurface $M_1\times M_2\times S^1$ immersed by $f$ is given by 
$$C(p_1,p_2,\iota(\theta))=-\frac{u_1'(\theta)^2-u_2'(\theta)^2}{u_1'(\theta)^2+u_2'(\theta)^2}.$$
From $(3.14)$, we can derive 
$$\begin{array}{l}
\hspace{0.5truecm}\displaystyle{df_{(p_1,p_2,\iota(\theta))}(A_{(p_1,p_2,\iota(\theta))}(\vv^1_{jk}))=-\widetilde{\nabla}^f_{\vv^1_{jk}}\overline{\vN}
=\frac{u_2'(\theta)}{\sqrt{u_1'(\theta)^2+u_2'(\theta)^2}}\,(\vY_1'(u_1(\theta)),{\bf 0})}\\
\displaystyle{=-\frac{u_2'(\theta)}{\sqrt{u_1'(\theta)^2+u_2'(\theta)^2}}\,\left(\sqrt{\mu_{1k,p_1}}\,\sin(u_1(\theta)\sqrt{\mu_{1k,p_1}})
+\lambda_{1j,p_1}\cos(u_1(\theta)\sqrt{\mu_{1k,p_1}})\right)}\\
\hspace{0.5truecm}\displaystyle{\left(P_{\gamma_{p_1}|_{[0,u_1(\theta)]}}((df_1)_{p_1}(\vv^1_{jk}),\,{\bf 0}\right)}
\end{array}\leqno{(3.15)}$$
and 
$$\begin{array}{l}
\hspace{0.5truecm}\displaystyle{df_{(p_1,p_2,\iota(\theta))}(A_{(p_1,p_2,\iota(\theta))}(\vv^2_{jk}))=-\widetilde{\nabla}^f_{\vv^2_{jk}}\overline{\vN}
=-\frac{u_1'(\theta)}{\sqrt{u_1'(\theta)^2+u_2'(\theta)^2}}\,({\bf 0},\,\vY_2'(u_2(\theta)))}\\
\displaystyle{=\frac{u_1'(\theta)}{\sqrt{u_1'(\theta)^2+u_2'(\theta)^2}}\,\left({\bf 0},\,\sqrt{\mu_{2k,p_2}}\,\sin(u_2(\theta)\sqrt{\mu_{2k,p_2}})
+\lambda_{2j,p_2}\cos(u_2(\theta)\sqrt{\mu_{2k,p_2}})\right)}\\
\hspace{0.5truecm}\displaystyle{\left({\bf 0},\,P_{\gamma_{p_2}|_{[0,u_2(\theta)]}}((df_2)_{p_2}(\vv^2_{jk})\right)}
\end{array}\leqno{(3.16)}$$
From $(3.8)$ and $(3.15)$, we obtain 
$$\begin{array}{l}
\hspace{0.5truecm}\displaystyle{df_{(p_1,p_2,\iota(\theta))}(A_{(p_1,p_2,\iota(\theta))}(\vv^1_{jk}))}\\
\displaystyle{=-\frac{u'_2(\theta)}{\sqrt{u'_1(\theta)^2+u'_2(\theta)^2}}\cdot
\frac{\sqrt{\mu_{1k,p_1}}\sin(u_1(\theta)\sqrt{\mu_{1k,p_1}})+\lambda_{1j,p_1}\cos(u_1(\theta)\sqrt{\mu_{1k,p_1}})}
{\cos(u_1(\theta)\sqrt{\mu_{1k,p_1}})-\frac{\lambda_{1j,p_1}\sin(u_1(\theta)\sqrt{\mu_{1k,p_1}})}{\sqrt{\mu_{1k,p_1}}}}\,
df_{(p_1,p_2,\iota(\theta))}(\vv^1_{jk})}
\end{array}$$
and hence 
$$\begin{array}{l}
\hspace{0.5truecm}\displaystyle{A_{(p_1,p_2,\iota(\theta))}(\vv^1_{jk})}\\
\displaystyle{=-\frac{u'_2(\theta)}{\sqrt{u'_1(\theta)^2+u'_2(\theta)^2}}\cdot
\frac{\sqrt{\mu_{1k,p_1}}\sin(u_1(\theta)\sqrt{\mu_{1k,p_1}})+\lambda_{1j,p_1}\cos(u_1(\theta)\sqrt{\mu_{1k,p_1}})}
{\cos(u_1(\theta)\sqrt{\mu_{1k,p_1}})-\frac{\lambda_{1j,p_1}\sin(u_1(\theta)\sqrt{\mu_{1k,p_1}})}{\sqrt{\mu_{1k,p_1}}}}\,
\vv^1_{jk}.}
\end{array}\leqno{(3.17)}$$
Similarly, from $(3.9)$ and $(3.16)$, we obtain 
$$\begin{array}{l}
\hspace{0.5truecm}\displaystyle{A_{(p_1,p_2,\iota(\theta))}(\vv^2_{jk})}\\
\displaystyle{=\frac{u'_1(\theta)}{\sqrt{u'_1(\theta)^2+u'_2(\theta)^2}}\cdot
\frac{\sqrt{\mu_{2k,p_2}}\sin(u_2(\theta)\sqrt{\mu_{2k,p_2}})+\lambda_{2j,p_2}\cos(u_2(\theta)\sqrt{\mu_{2k,p_2}})}
{\cos(u_2(\theta)\sqrt{\mu_{2k,p_2}})-\frac{\lambda_{2j,p_2}\sin(u_2(\theta)\sqrt{\mu_{2k,p_2}})}{\sqrt{\mu_{2k,p_2}}}}\,
\vv^2_{jk}.}
\end{array}\leqno{(3.18)}$$
From $(3.5)$ and $(3.14)$, we can derive 
$$\begin{array}{l}
\hspace{0.5truecm}\displaystyle{df_{(p_1,p_2,\iota(\theta))}\left(A_{(p_1,p_2,\iota(\theta))}\left(\frac{\partial}{\partial\theta}\right)\right)
=-\widetilde{\nabla}^f_{\frac{\partial}{\partial\theta}}\overline{\vN}}\\
\displaystyle{=\left(\left(\frac{-u'_2(\theta)}{\sqrt{u'_1(\theta)^2+u'_2(\theta)^2}}\right)'\,P_{\gamma_{p_1}|_{[0,u_1(\theta)]}}((\vN_1)_{p_1}),\right.}\\
\hspace{0.9truecm}\displaystyle{\left.\left(\frac{u'_1(\theta)}{\sqrt{u'_1(\theta)^2+u'_2(\theta)^2}}\right)'\,
P_{\gamma_{p_2}|_{[0,u_2(\theta)]}}((\vN_2)_{p_2})\right)}\\
\displaystyle{=\frac{u'_1(\theta)u''_2(\theta)-u''_1(\theta)u'_2(\theta)}{(u'_1(\theta)^2+u'_2(\theta)^2)^{\frac{3}{2}}}\,
df_{(p_1,p_2,\iota(\theta))}\left(\frac{\partial}{\partial\theta}\right)}
\end{array}$$
and hence 
$$A_{(p_1,p_2,\iota(\theta))}\left(\frac{\partial}{\partial\theta}\right)
=\frac{u'_1(\theta)u''_2(\theta)-u''_1(\theta)u'_2(\theta)}{(u'_1(\theta)^2+u'_2(\theta)^2)^{\frac{3}{2}}}\,
\frac{\partial}{\partial\theta}.\leqno{(3.19)}$$
From $(3.14)$, we have 
$$\begin{array}{l}
\displaystyle{\widetilde R(\overline{\vN})_{(p_1,p_2,\iota(\theta))}(\vv^1_{jk})=\frac{u'_2(\theta)^2\mu_{1k,p_1}}{u'_1(\theta)^2+u'_2(\theta)^2}\,\vv^1_{jk},}\\
\displaystyle{\widetilde R(\overline{\vN})_{(p_1,p_2,\iota(\theta))}(\vv^2_{jk})=\frac{u'_1(\theta)^2\mu_{2k,p_2}}{u'_1(\theta)^2+u'_2(\theta)^2}\,\vv^2_{jk}.}
\end{array}\leqno{(3.20)}$$
Also, from $(3.5)$ and $(3.14)$, we have 
$$\widetilde R(\overline{\vN})_{(p_1,p_2,\iota(\theta))}\left(\frac{\partial}{\partial\theta}\right)={\bf 0}.\leqno{(3.21)}$$
From $(3.17),\,(3.18),\,(3.19),(3.20)$ and $(3.21)$, the common eigenspace decomposition of $A_{(p_1,p_2,\iota(\theta))}$ and $\widetilde R(\overline{\vN})_{(p_1,p_2,\iota(\theta))}$ 
is given by 
$$\begin{array}{l}
\hspace{0.5truecm}\displaystyle{T_{(p_1,p_2,\iota(\theta))}(M_1\times M_2\times S^1)=T_{p_1}M_1\oplus T_{p_2}M_2\oplus{\rm Span}\left\{\frac{\partial}{\partial\theta}\right\}}\\
\displaystyle{=\left(\mathop{\oplus}_{(j,k)\in I_{A_1R_1,p_1}}E^1_{jk,p_1}\right)\oplus\left(\mathop{\oplus}_{(j,k)\in I_{A_2R_2,p_2}}E^2_{jk,p_2}\right)
\oplus{\rm Span}\left\{\frac{\partial}{\partial\theta}\right\}.}
\end{array}\leqno{(3.22)}$$
Hence the hypersurface $M_1\times M_2\times S^1$ immersed by $f$ is curvature-adapted.  \qed

%\vspace{0.5truecm}
%\centerline{\input{CAMFS-F1.tex}\hspace{0truecm}}
%\vspace{0.3truecm}
%\centerline{\bf Figure 1.}
%\vspace{0.5truecm}

\vspace{0.5truecm}

{\small 
\rightline{Department of Mathematics, Faculty of Science}
\rightline{Tokyo University of Science, 1-3 Kagurazaka}
\rightline{Shinjuku-ku, Tokyo 162-8601 Japan}
\rightline{(koike@rs.tus.ac.jp)}
}

\end{document}